\documentclass{article}
\usepackage{graphicx} % Required for inserting images
\usepackage{amsmath}
\usepackage{amssymb}
\usepackage{subcaption}
\usepackage{booktabs}
\usepackage{authblk}

\title{Polynomial Optimization via Random Projection and Consensus}
\author{Etienne Buehrle}
\author{Christoph Stiller}
\affil{Karlsruhe Institute of Technology\\\texttt{\{etienne.buehrle, stiller\}@kit.edu}}
\date{}

\def\A{\mathcal A}
\def\X{\mathcal X}
\def\U{\mathcal U}
\def\S{\mathbb S}
\def\R{\mathbb R}
\DeclareMathOperator{\Tr}{Tr}

\begin{document}

\maketitle

\begin{abstract}
    We propose a black-box approach to reducing large semidefinite programs to a set of smaller semidefinite programs by projecting to random linear subspaces. We evaluate our method on a set of polynomial optimization problems, demonstrating improved scalability.
\end{abstract}

\section{Introduction}
Polynomial programs arise in a number of applications, including control theory and combinatorial optimization \cite{parrilo2012chapter}. They can be efficiently represented as semidefinite programs, which are themselves a common class of convex optimization problems and can be solved in polynomial time using off-the-shelf interior-point solvers \cite{boyd2004convex}. However, in semidefinite optimization, the computational complexity is cubic in the number of variables, making their practical use difficult beyond problems of modest size.

In certain cases, polynomial programs can be reduced by leveraging symmetry \cite{gatermann2004symmetry} and sparsity properties \cite{magron2023sparse} of the polynomial variables. These approaches, however, require knowledge of the objective and constraints and do not apply to all objective functions and constraints.

Prior works have also proposed to use inner approximations to the positive semidefinite cone \cite{ahmadi2019dsos}, which restrict the space of solutions to computationally more advantageous cones. While alternative approaches including natural cone formulations of the moment cone \cite{papp2019sum, coey2022solving} are able to avoid this restriction, they generally do not allow direct access to the polynomial variable and do not permit specifying constraints on its coefficients. In contrast, the semidefinite formulation deals with polynomials in an explicit representation, allowing for easy manipulations.

We propose an alternative approach that keeps the semidefinite structure largely intact, but projects the positive semidefiniteness constraint onto a set of randomly sampled linear subspaces, reducing the problem dimensionality. The optimization variable is reconstructed by enforcing consistency constraints between the projections, akin to compressed sensing \cite{candes2005decoding}.

\subsection{Polynomial Optimization}
Polynomial optimization problems can be cast as semidefinite programs by strengthening the polynomial positivity constraints to sums of squares constraints \cite{parrilo2012chapter}. A typical polynomial optimization problem is the minimization of a polynomial $p : \X \to \R$ \begin{equation}\label{eq:pop}
    \begin{split}
        \underset{\lambda \in \R}{\text{maximize}} \quad& \lambda \\
        \text{subject to} \quad& p-\lambda \in SOS(\X)
    \end{split}
\end{equation} which can be interpreted as maximizing a lower bound $\lambda$ to the polynomial.

The sum of squares constraint $p \in SOS(\X)$ amounts to finding a positive semidefinite Gram matrix $Q \succeq 0$ such that $p(x) = \phi(x)^\top Q \phi(x)$ for a basis\footnote{We use the canonical basis $\phi(x) = \begin{bmatrix}1 & x & \dots & x^d\end{bmatrix}^\top$ of bounded degree monomials.} $\phi$. This is sufficient to ensure that $p(x) \geq 0$ for all $x \in \X$, leading to a semidefinite program.

Although convex, semidefinite programs are among the most challenging class of optimization problems, with computational complexity cubic in the dimensionality of the positive semidefinite cone \cite{boyd2004convex}. In the following, we propose a dimensionality reduction method based on projections onto random linear subspaces.

\section{Background}
We consider semidefinite programs of the form \begin{equation}\label{eq:sdp}\begin{split}
    \underset{X \in \S^n}{\text{minimize}} \quad& \Tr(CX) \\
    \text{subject to} \quad& \A(X) = b \\
    & X \succeq 0
\end{split}\end{equation} with symmetric $C \in \R^{n \times n}$, linear map $\A : \S^n \to \R^m$, vector $b \in \R^m$, and $\S^n$ the set of symmetric $n \times n$ matrices. The computational cost to solve such problems using interior-point methods is cubic in the dimension $n$ of the positive semidefinite cone \cite{boyd2004convex}.

\subsection{Dual Semidefinite Program}
The dual problem to \eqref{eq:sdp} is \begin{equation}\label{eq:sdp-dual}
    \begin{split}
        \underset{y\in\R^m, S \succeq 0}{\text{maximize}} \quad& b^\top y \\
        \text{subject to} \quad& \A^*(y) + S = C
    \end{split}
\end{equation} with dual variables $y$ and $S \in \S^n$. Since the semidefinite cone is self-dual, the dual problem is a semidefinite program.

\section{Method}
We modify \eqref{eq:sdp} by picking $N$ random linear subspaces $U_i \in \R^{n \times r}$, $i=1, \dots, N$, $r < n$, and note that $U^\top X U \succeq 0$ if $X \succeq 0$, yielding \begin{equation}\label{eq:csdp}\begin{split}
    \underset{X \in \S^n}{\text{minimize}} \quad& \Tr(CX) \\
    \text{subject to} \quad& \A(X) = b \\
    & U_i^\top X U_i \succeq 0 \quad i=1, \dots, N
\end{split}\end{equation} where positive semidefiniteness is now imposed on the smaller cone of size $r \times r$. The computational cost is cubic in the rank $r$ of the projections, at the expense of solving $N$ optimization problems, which can be parallelized \cite{boyd2011distributed}.

For the reverse implication, note that if there exists $v$ such that $v^\top Xv < 0$, then there exist $U$, $\hat{v}$ such that $U\hat{v}=v$ and $\hat{v}^\top U^\top X U \hat{v} < 0$. Such a $U$ is increasingly likely to be sampled as $r$ increases.

\subsection{Dual Semidefinite Program}
The dual problem to \eqref{eq:csdp} can be written as the semidefinite program \begin{equation}\label{eq:csdp-dual}
    \begin{split}
        \underset{y\in\R^m, S_i \succeq 0}{\text{maximize}} \quad& b^\top y \\
        \text{subject to} \quad& \A^*(y) + \sum_{i=1}^N U_iS_iU_i^\top = C
    \end{split}
\end{equation} with duals $y$ and $S_i\in\S^r$, $i=1,\dots,N$, where the dual variable $S$ in \eqref{eq:csdp} is now approximated by the sum of outer products $\sum_{i=1}^N U_iS_iU_i^\top$.

\section{Results}
We evaluate our method on a synthetic polynomial optimization problem and a real-world polynomial optimal control problem.

\subsection{Polynomial Optimization}
\begin{figure}
    \centering
    \begin{subfigure}{0.25\linewidth}
        \includegraphics[width=\linewidth]{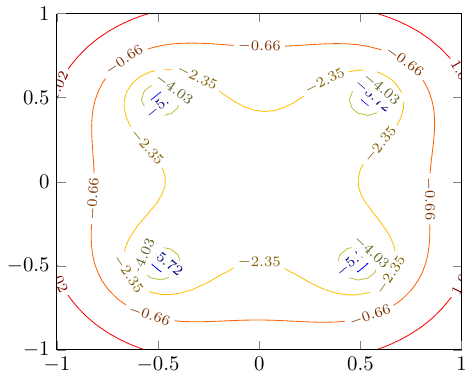}
        \caption{Polynomial}
        \label{fig:pop-p}
    \end{subfigure}
    \hfill
    \begin{subfigure}{0.3\linewidth}
        \includegraphics[width=\linewidth]{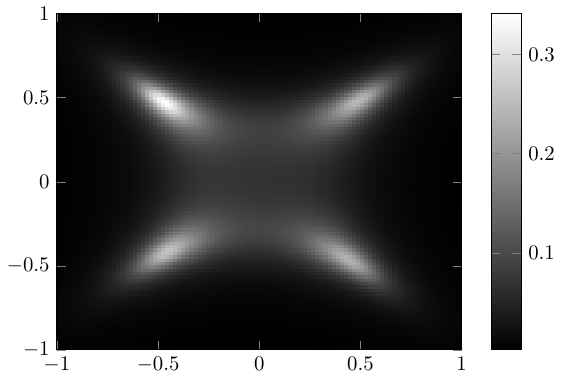}
        \caption{$r = 6$}
    \end{subfigure}
    \hfill
    \begin{subfigure}{0.325\linewidth}
        \includegraphics[width=\linewidth]{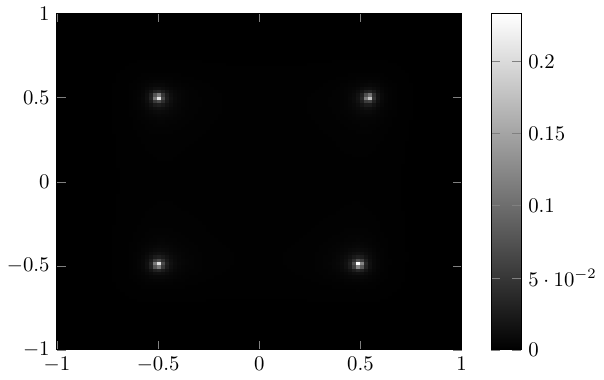}
        \caption{$r = 11$}
    \end{subfigure}
    \caption{Polynomial optimization with varying ranks ($N=100$).}
    \label{fig:pop}
\end{figure}

We consider a polynomial optimization problem of the form \eqref{eq:pop} on a polynomial of degree 8 in two variables with four double zeros (Figure \ref{fig:pop-p}), leading to a sum of squares constraint in 495 coefficients.

\begin{table}
    \centering
    \begin{tabular}{l|c|c|c|c|c|c|c|c|c|c}
        \toprule
        Rank $r$  & 1 & 3 & 6 & 9 & 11 & 14 & 17 & 19 & 22 & 25 \\
        \midrule
        Cone size & 1 & 6 & 21 & 45 & 66 & 105 & 153 & 190 & 253 & 325 \\
        Objective & $-\infty$ & $-0.3$ & 0 & 0 & 0 & 0 & 0 & 0 & 0 & 0 \\
        \bottomrule
    \end{tabular}
    \caption{Cone sizes and attained objective values for varying ranks.}
    \label{tab:pop}
\end{table}

We fix a sample size $N=100$ and show the attained objective value for different projection ranks $r$ in Table \ref{tab:pop}, finding that the objective value remains stable for wide ranges of $r$. The size of the positive semidefinite cone is determined by the rank $r$ of the projections, effectively reducing the problem dimensionality.

We show the dual of the polynomial program in Figure \ref{fig:pop}, which can be interpreted as a probability density over the minimizer (cf. Appendix \ref{sec:pop-dual}).

\subsection{Polynomial Optimal Control}
We consider a polynomial optimal control problem of a dynamical system $\dot x = f(x,u)$, $x(t) \in \X$, $u(t) \in \U$, subject to a running cost $c: \X \times \U \to \R$, which can be formulated as a subsolution of the Bellman equation \begin{equation}\label{eq:poc}\begin{split}
    \underset{V \in SOS(\X)}{\text{maximize}} \quad& V(x_0) - V(x_T) \\
    \text{subject to} \quad& \nabla V^\top f + c \in SOS(\X \times \U)
\end{split}\end{equation} with initial state $x(0) = x_0$, time horizon $T$, and terminal state $x(T) = x_T$. The polynomial variable $V$ can be interpreted as a value function \cite{lasserre2008nonlinear}.

\begin{figure}
    \centering
    \begin{subfigure}{0.3\linewidth}
        \centering
        \includegraphics[width=\linewidth]{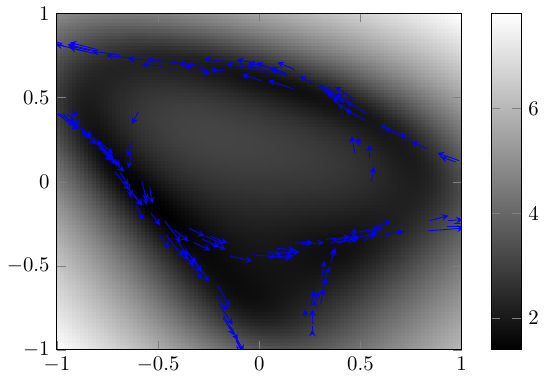}
        \caption{Marginal cost}
        \label{fig:poc-p}
    \end{subfigure}
    \hfill
    \begin{subfigure}{0.3\linewidth}
        \centering
        \includegraphics[width=\linewidth]{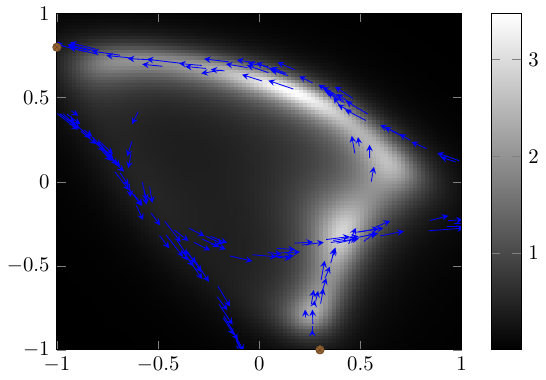}
        \caption{$r=14$}
    \end{subfigure}
    \hfill
    \begin{subfigure}{0.3\linewidth}
        \centering
        \includegraphics[width=\linewidth]{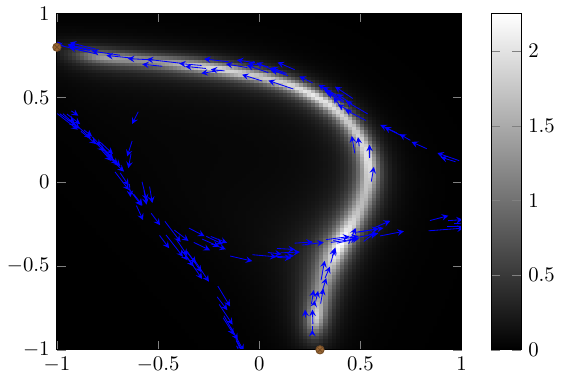}
        \caption{$r=25$}
    \end{subfigure}
    \caption{Polynomial optimal control with varying ranks ($N=100$).}
    \label{fig:poc}
\end{figure}

\begin{table}
    \centering
    \begin{tabular}{l|c|c|c|c|c|c|c|c|c}
        \toprule
        Rank $r$  & 1 & 3 & 6 & 9 & 11 & 14 & 17 & 19 & 25 \\
        \midrule
        Objective & $-\infty$ & $-\infty$ & $-\infty$ & $-\infty$ & $-\infty$ & 151.9 & 177.0 & 177.6 & 177.6 \\
        \bottomrule
    \end{tabular}
    \caption{Attained objective value for varying ranks.}
    \label{tab:poc}
\end{table}

We choose a cost function representing a segment of a traffic roundabout \cite{zhan2019interaction, buehrle2025stochastic}, as shown in Figure \ref{fig:poc-p}. We fix a sample size $N=100$ and show the attained objective value for varying projection ranks $r$ in Table \ref{tab:poc}, finding that, as in the case of polynomial optimization, the objective value remains stable for wide ranges of $r$.

We show the dual of the polynomial program in Figure \ref{fig:poc}, which can be interpreted as an occupation measure of the state trajectory (cf. Appendix \ref{sec:poc-dual}).

\section{Conclusion}
We have presented an approach to reducing large semidefinite programs to a set of smaller semidefinite programs. Our approach is black-box, requiring no knowledge about the internal problem structure, while keeping the modeling and solution process largely intact. In evaluations on synthetic and real-world example problems, we find that robust approximations can be obtained using small projection ranks, leading to effective reductions in the dimensionality and computational complexity. 

Future work will explore automatic ways to choose the projection rank, as well as improved subspace sampling and dual recovery methods.

\bibliographystyle{plain}
\bibliography{main}

\appendix
\def\M{\mathcal M}

\section{Dual Polynomial Program}\label{sec:pop-dual}
The dual of the polynomial program \eqref{eq:pop} is a measure program \begin{equation}\label{eq:pop-dual}\begin{split}
    \underset{\mu\in\M_+(\X)}{\text{minimize}} \quad& \langle p, \mu \rangle \\
    \text{subject to} \quad& \langle 1, \mu \rangle = 1
\end{split}\end{equation} with inner product $\langle p, \mu \rangle = \int_\X p\,d\mu$, which can be interpreted as the continuous limit of a weighted sum \cite{lasserre2009moments}.

\section{Dual Polynomial Optimal Control}\label{sec:poc-dual}
The dual of the polynomial optimal control problem \eqref{eq:poc} is a measure program \begin{equation}\begin{split}
    \underset{\mu\in\M_+(\X \times \U)}{\text{minimize}} \quad& \langle c, \mu \rangle \\
    \text{subject to} \quad& \nabla\cdot(f\mu) = \delta_{x_0} - \delta_{x_T}
\end{split}\end{equation} where $\delta_x$ denotes the Dirac measure at $x$. The measure $\mu$ can be interpreted as the occupation measure of a state-action trajectory subject to a continuity equation, with a source at the initial state and a sink at the terminal state \cite{lasserre2009moments}.

\end{document}